\documentclass[leqno]{article}
\usepackage{geometry}
\usepackage{graphicx}	
\usepackage{tikz}
\usepackage[utf8]{inputenc}
\usepackage[10pt]{type1ec}
\usepackage[T1,T2A]{fontenc}
\usepackage[protrusion=true,expansion=true]{microtype}
\usepackage[russian,english]{babel}
\usepackage{mathtools}
\usepackage{amsfonts,amssymb,amsmath,amsthm}

\usepackage[colorinlistoftodos]{todonotes}

\usepackage{url}
\usepackage{parskip}
\ifpdf
\usepackage{hyperref}
\else
\usepackage[dvipdfmx]{hyperref}
\fi

\usepackage{stmaryrd}

\usepackage{csquotes}
\usepackage[backend=biber,style=alphabetic,autolang=other,giveninits=true,maxbibnames=99]{biblatex}

\DeclareNameAlias{default}{family-given}
\usepackage{etoolbox}
\pdfpageheight=\textheight
\pdfpagewidth=\textwidth
\hypersetup{colorlinks=true}
\def\thtext#1{
  \catcode`@=11
  \gdef\@thmcountersep{. #1}
  \catcode`@=12
}

\def\threst{
  \catcode`@=11
  \gdef\@thmcountersep{.}
  \catcode`@=12
}

\theoremstyle{plain}

\newtheorem{theorem}{Theorem}[section]

\newtheorem{lemma}[theorem]{Lemma}

\theoremstyle{definition}
\newtheorem{define}[theorem]{Definition}

 \catcode`@=11
 \def\.{.\spacefactor\@m}
 \catcode`@=12

\def\Z{{\mathbb Z}}

\let\textcaron\v 
\AtBeginBibliography{\let\v\textcaron}
\def\v{\varphi}

\def\0{\emptyset}
\def\:{\colon}
\def\<{\langle}
\def\>{\rangle}
\def\[{\llbracket}
\def\]{\rrbracket}

\def\rom#1{\emph{#1}}
\def\({\rom(}
\def\){\rom)}

\def\RR{\mathcal{R}}

\newcommand{\dGH}[1]{\operatorname{d}^{#1}_{\mathrm{GH}}}

\def\dis{\operatorname{dis}}

\def\GH{\operatorname{\mathcal{G\!H}}}

\def\gr{\operatorname{gr}}

\def\id{\operatorname{id}}

\begin{document}
\title{Bounded Spaces without a Geodesic\\ in the Gromov--Hausdorff Class\thanks{The work was supported by the Basis Foundation (grant No.~25-8-3-11-1).}}
\author{A.~A.~Vikhrov}
\date{}
\maketitle
\begin{abstract}
We construct two countable complete bounded metric spaces at Gromov--Hausdorff distance $1/4$ with no midpoint and hence no shortest curve between them.
The proof uses the local existence of optimal correspondences near a space with finitely many self-correspondences of distortion below a fixed positive threshold, each of which is the graph of a bijection.
\end{abstract}

\section{Introduction}\label{sec:history}
The Gromov--Hausdorff distance between nonempty metric spaces $A$ and $B$ is the infimum of the Hausdorff distances between their isometric copies in a common metric space:
$$
 \dGH{}(A,B)=\inf_{Z,f,g}d_H^Z\bigl(f(A),g(B)\bigr),
$$
where $f:A\to Z$ and $g:B\to Z$ are isometric embeddings.
For nonempty subsets $U,V\subset Z$, the Hausdorff distance $d_H^Z(U,V)$ is the least $r\ge0$ such that every point of $U$ is at distance at most $r$ from $V$, and every point of $V$ is at distance at most $r$ from $U$; if there is no such $r$, we set $d_H^Z(U,V)=\infty$.

We consider the class $\GH$ of all metric spaces up to isometry, equipped with the Gromov--Hausdorff distance.
A shortest curve, or geodesic, is a curve whose finite length equals the distance between its endpoints.
We study the existence of shortest curves between arbitrary metric spaces at finite Gromov--Hausdorff distance.

For compact metric spaces, the existence of geodesics was proved by Ivanov, Nikolaeva and Tuzhilin: they constructed a midpoint for each pair and used the completeness of the space of compact metric spaces~\cite[\textsection~4]{INT}.

Another proof gave an explicit construction of a geodesic.
Ivanov, Iliadis and Tuzhilin, and independently Chowdhury and Mémoli, proved that any two nonempty compact metric spaces admit an optimal correspondence~\cites[Theorem~2.6]{IIT}[Proposition~1.1]{CM}.
They then used such a correspondence to construct a geodesic~\cites[Theorem~3.1 and Introduction]{IIT}[Theorem~1.2]{CM}.

Let $R$ be an optimal correspondence between $A$ and $B$.
For $0<t<1$, define a metric on $R$ by
\begin{equation}\label{eq:linear}
 d_t\bigl((a,b),(a',b')\bigr)=(1-t)d_A(a,a')+t\,d_B(b,b').
\end{equation}
The curve $\gamma$ with $\gamma(0)=A$, $\gamma(1)=B$ and $\gamma(t)=(R,d_t)$ for $0<t<1$ is called a linear geodesic.

This construction also applies to noncompact spaces whenever they admit an optimal correspondence~\cite[Theorem~2.1 with $\varepsilon=0$]{BIT}.
Moreover, Borzov, Ivanov and Tuzhilin proved that, for spaces $A,B$ at finite Gromov--Hausdorff distance and every $\varepsilon>0$, using a correspondence of distortion at most $2\dGH{}(A,B)+2\varepsilon$ in the same formula gives a curve of length at most $\dGH{}(A,B)+\varepsilon$~\cite[Theorem~2.1]{BIT}.
Thus any two spaces at finite distance can be joined by curves whose lengths are arbitrarily close to the distance between them.

Vikhrov constructed a pair consisting of a two-point space of diameter $1$ and a countable complete metric space with no optimal correspondence between them~\cite[Example~4.1, Lemma~4.5]{Vikhrov2026}.
For such a pair, a linear geodesic cannot be constructed by formula~\eqref{eq:linear}.
Nor can a shortest curve be obtained by joining finitely many linear geodesics whose adjacent endpoints are isometric: the existence of such a piecewise linear geodesic implies the existence of an optimal correspondence between its endpoints~\cite[Proposition~3.2]{Vikhrov2026}.
These are the spaces denoted by $A,B$ in~\cite[Algorithm~1]{Vikhrov2025} and by $X,Y$ in Example~4.1 of the 2026 paper.
They also have no optimal Hausdorff realization: the two spaces cannot be isometrically embedded in a common metric space so that the Hausdorff distance between their images equals their Gromov--Hausdorff distance~\cite[Theorem~6]{Vikhrov2025}.

Nevertheless, a geodesic between these spaces was constructed in~\cite[Theorem~5.6]{Vikhrov2026}.
Its construction uses two linear geodesics whose adjacent endpoints are nonisometric spaces at zero Gromov--Hausdorff distance.
Thus the absence of an optimal correspondence and an optimal Hausdorff realization does not rule out a shortest curve.

In this paper we construct countable complete bounded spaces $X,Y$ with no shortest curve between them: they have no midpoint (Theorem~\ref{thm:bounded}).
They admit no optimal correspondence, but, unlike the earlier example, the existence of a midpoint would imply the existence of such a correspondence.
Indeed, Lemma~\ref{lem:local} gives optimal correspondences from both endpoints to the supposed midpoint; composing one with the inverse of the other gives an optimal correspondence between the endpoints.
For the constructed pair,
$$
 \dGH{}(X,Y)=\frac14,\qquad
 \nexists M\in\GH:\quad \dGH{}(X,M)=\dGH{}(M,Y)=\frac18.
$$

AI tools were used during the preparation of this article.

\section{Preliminaries}\label{sec:prelim}
We introduce the notation used in the proof of Theorem~\ref{thm:bounded}.
All metric spaces considered below are nonempty.

\begin{define}\label{def:class}
The Gromov--Hausdorff class $\GH$ is the proper class, in von Neumann--Bernays--Gödel set theory, of all metric spaces considered up to isometry.
\end{define}
The Gromov--Hausdorff distance is a generalized pseudometric on $\GH$: infinite values and zero distance between nonisometric spaces are allowed~\cite[Corollary~2.2]{BIT}.

\begin{define}\label{def:corr}
A correspondence between $A$ and $B$ is a relation $R\subset A\times B$ whose two projections are surjective.
We write $\RR(A,B)$ for the set of all such correspondences.
The distortion of a correspondence $R$ is the quantity
$$
 \dis R=\sup\bigl\{|d_A(a,a')-d_B(b,b')|:(a,b),(a',b')\in R\bigr\}.
$$
\end{define}
For any metric spaces $A,B$, we have (see Burago, Burago and Ivanov~\cite[Theorem~7.3.25]{BBI})
\begin{equation}\label{eq:gh}
 2\dGH{}(A,B)=\inf_{R\in\RR(A,B)}\dis R.
\end{equation}
A correspondence $R\in\RR(A,B)$ is called optimal if $\dis R=2\dGH{}(A,B)<\infty$.

For $R\in\RR(A,B)$ and $T\in\RR(B,C)$, set
$$
 T\circ R=\{(a,c):(a,b)\in R\text{ and }(b,c)\in T\text{ for some }b\in B\}.
$$
Then $T\circ R\in\RR(A,C)$ and
\begin{equation}\label{eq:composition}
 \dis(T\circ R)\le\dis T+\dis R.
\end{equation}
The inverse correspondence $R^{-1}=\{(b,a):(a,b)\in R\}\in\RR(B,A)$ has the same distortion as $R$.

\begin{define}\label{def:geodesic}
The length of a continuous curve $\gamma:[0,1]\to\GH$ is the quantity
$$
 L(\gamma)=\sup_{0=t_0<\dots<t_N=1}
 \sum_{i=1}^N\dGH{}\bigl(\gamma(t_{i-1}),\gamma(t_i)\bigr).
$$
A curve $\gamma$ is called a geodesic, or a shortest curve, if $L(\gamma)=\dGH{}(\gamma(0),\gamma(1))<\infty$.
A midpoint between $A$ and $B$ at finite distance is a space $M$ such that $\dGH{}(A,M)=\dGH{}(M,B)=\dGH{}(A,B)/2$.
\end{define}
Every shortest curve passes through a midpoint between its endpoints.

\section{Local Existence of Optimal Correspondences}\label{sec:local}
We prove a lemma that will rule out a midpoint between the spaces of the counterexample.
For a metric space $A$ and $c>0$, write
$$
 \RR_{<c}(A,A)=\{C\in\RR(A,A):\dis C<c\}.
$$

\begin{lemma}\label{lem:local}
Let $A$ be a metric space and $c>0$, and suppose that $\RR_{<c}(A,A)$ is finite and consists of graphs of bijections.
Then every metric space $M$ with $r=\dGH{}(A,M)<c/4$ admits an optimal correspondence $R\in\RR(A,M)$.
\end{lemma}
\begin{proof}
For a map $g:A\to A$, write $\gr g=\{(a,g(a)):a\in A\}$ for its graph.
Write $\RR_{<c}(A,A)=\{\gr g_1,\dots,\gr g_N\}$, where $g_1,\dots,g_N:A\to A$ are bijections.
Choose $\varepsilon_n>0$ such that $\varepsilon_n\to0$ and $4r+2\sup_n\varepsilon_n<c$.
For each $n$, choose $R_n\in\RR(A,M)$ with $\dis R_n<2r+\varepsilon_n$.
By the composition estimate~\eqref{eq:composition}, the self-correspondences $C_n=R_n^{-1}\circ R_1$ satisfy
$$
 \dis C_n\le\dis R_n+\dis R_1
 <4r+\varepsilon_n+\varepsilon_1<c.
$$
Thus $C_n\in\RR_{<c}(A,A)$, so $C_n=\gr g_{i_n}$ for some $i_n\in\{1,\dots,N\}$.
For each $a\in A$, there is $m\in M$ with $(a,m)\in R_1$, so $(a,a)\in R_1^{-1}\circ R_1=C_1$.
Since $C_1$ is the graph of a map and contains the entire diagonal, it equals the diagonal.

For each $m\in M$, there is at least one $a\in A$ with $(a,m)\in R_1$.
This point is unique: if $(a,m),(a',m)\in R_1$, then $(a,a')\in C_1$, whence $a=a'$.
Write $p(m)$ for this point.
Thus $R_1^{-1}$ is the graph of a map $p:M\to A$.

Now let $(a,m)\in R_n$.
Then
$$
 (p(m),m)\in R_1,\qquad (m,a)\in R_n^{-1},
$$
so $(p(m),a)\in R_n^{-1}\circ R_1=C_n$ and $a=g_{i_n}(p(m))$.
Since the projection of $R_n$ onto $M$ is surjective, such a pair exists for every $m$, and we get
$$
 R_n=\{(g_{i_n}(p(m)),m):m\in M\}.
$$
Since $p$ does not depend on $n$ and $i_n\in\{1,\dots,N\}$, there are at most $N$ distinct relations $R_n$.
Some $R\in\RR(A,M)$ occurs along an infinite subsequence $R_{n_k}$.
Then $\dis R<2r+\varepsilon_{n_k}$ for all $k$, hence $\dis R\le2r$.
The reverse inequality follows from the distance formula~\eqref{eq:gh}, so $R$ is optimal.
\end{proof}

\section{A Bounded Counterexample}\label{sec:counterexample}
We construct a pair with no midpoint.
Both spaces have underlying set $\Z$.
Write $n\wedge m=\min\{n,m\}$ and $e(j)=1/(|j|+2)$ for $j\in\Z$.
Define the distances by the following table, and set them equal to zero on the diagonal:
\begin{equation}\label{eq:metrics}
\begin{array}{c|c|c}
 &d_X(n,m)&d_Y(n,m)\\ \hline
 |n-m|=1&3+\mathbf1_{\{0\}}(n\wedge m)&7/2-e(n\wedge m)\\
 |n-m|=2&11/2&5+\mathbf1_{\{0\}}(n\wedge m)\\
 |n-m|\ge3&11/2&11/2
\end{array}
\end{equation}
Here $\mathbf1_{\{0\}}(j)$ equals $1$ if $j=0$ and $0$ otherwise.
The term $e(j)$ will allow the distortions of shifts to approach their infimum without attaining it.
The distinguished pairs $\{0,1\}$ in $X$ and $\{0,2\}$ in $Y$ will ensure that each space has only finitely many self-correspondences of distortion less than $1$.

All nonzero distances in~\eqref{eq:metrics} lie in $[3,6]$.
For three distinct points, the triangle inequality follows from
$$
 d_A(n,m)\le6\le3+3\le d_A(n,p)+d_A(p,m),\qquad A\in\{X,Y\}.
$$
If some points coincide, the inequality is obvious; symmetry and positivity follow from the formulas.
Since distinct points in either space are at distance at least $3$, every Cauchy sequence is eventually constant, so both spaces are complete.
They are countable and bounded.

\begin{lemma}\label{lem:classification}
Every correspondence $R\in\RR(A,B)$ of distortion less than $1$, where $A,B\in\{X,Y\}$, is the graph of a bijection of the form
\begin{equation}\label{eq:maps}
 \varphi(n)=n+k\quad\text{or}\quad\varphi(n)=-n+k,\qquad k\in\Z.
\end{equation}
\end{lemma}
\begin{proof}
Let $A,B\in\{X,Y\}$ and $R\in\RR(A,B)$ with $\dis R<1$.
If one point on either side is related by $R$ to two distinct points on the other, then $\dis R\ge3$.
Thus $R$ is the graph of a bijection $\varphi:\Z\to\Z$.
Distances between adjacent integers lie in $[3,4]$, whereas those between distinct nonadjacent integers lie in $[5,6]$.
Since $\dis R<1$, the bijection $\varphi$ preserves adjacency in both directions.
Since $\varphi(n-1)$ and $\varphi(n+1)$ are distinct neighbors of $\varphi(n)$, we have $\varphi(n+1)-\varphi(n)=\varphi(n)-\varphi(n-1)$ for every $n\in\Z$, hence $\varphi(n)=\varphi(0)+n(\varphi(1)-\varphi(0))$.
This gives the stated form of $\varphi$.
\end{proof}

\begin{theorem}\label{thm:bounded}
Formulas~\eqref{eq:metrics} define countable complete bounded metric spaces $X,Y$ of diameter at most $6$ such that $\dGH{}(X,Y)=1/4$ and there is no metric space $M$ with
$$
 \dGH{}(X,M)=\dGH{}(M,Y)=\frac18.
$$
In particular, there is no shortest curve between $X$ and $Y$ in $\GH$.
\end{theorem}

\begin{proof}
We first compute the distortions of the graphs of the bijections~\eqref{eq:maps} from $X$ to $Y$.
Write $R_k$ and $S_k$ for the graphs of the maps $n\mapsto n+k$ and $n\mapsto-n+k$, respectively.
Both maps preserve $|n-m|$.
For a pair $\{j,j+q\}$ with $q\ge1$, let $h$ be the smaller endpoint of its image; then
$$
 h=\begin{cases}j+k,&\text{for }R_k,\\ k-j-q,&\text{for }S_k.\end{cases}
$$
If $q\ge3$, both distances equal $11/2$, so their difference is zero.
If $q=2$, the absolute difference is
$$
 \left|\frac{11}2-\left(5+\mathbf1_{\{0\}}(h)\right)\right|=\frac12.
$$
For an adjacent pair $\{j,j+1\}$ with $j\ne0$,
$$
 \left|d_X(j,j+1)-d_Y(h,h+1)\right|
 =\left|3-\left(\frac72-e(h)\right)\right|
 =\frac12-e(h)\le\frac12,
$$
since $0<e(h)\le1/2$.
Thus only the pair $\{0,1\}$ can give a difference larger than $1/2$.
Its image has smaller endpoint $k$ for $R_k$ and $k-1$ for $S_k$, so
$$
 \begin{aligned}
 \dis R_k&=\left|4-\left(\frac72-e(k)\right)\right|=\frac12+e(k),\\
 \dis S_k&=\left|4-\left(\frac72-e(k-1)\right)\right|=\frac12+e(k-1).
 \end{aligned}
$$
Both distortions are greater than $1/2$.
By Lemma~\ref{lem:classification}, all other correspondences have distortion at least $1$.
As $k$ ranges over $\Z$, so does $k-1$.
Since $e(k)>0$ and $e(k)\to0$ as $|k|\to\infty$, the distortions in both families have infimum $1/2$.
By the distance formula~\eqref{eq:gh},
\begin{equation}\label{eq:distance}
 2\dGH{}(X,Y)=\inf_{k\in\Z}\left(\frac12+e(k)\right)=\frac12,
\end{equation}
and no correspondence has distortion at most $1/2$.

We now verify the hypotheses on $A$ in Lemma~\ref{lem:local} with $c=1$, for $A=X$ and $A=Y$.
Let $R\in\RR(A,A)$, where $A\in\{X,Y\}$, satisfy $\dis R<1$.
By Lemma~\ref{lem:classification}, $R=\gr\varphi$, where $\varphi$ has one of the forms~\eqref{eq:maps}.
For $A=X$, we have $d_X(0,1)=4$, whereas $d_X(j,j+1)=3$ for $j\ne0$.
Thus the bijection $\varphi$ in~\eqref{eq:maps} must preserve the pair $\{0,1\}$: otherwise its image is another adjacent pair, and
$$
 \dis R\ge|d_X(0,1)-d_X(\varphi(0),\varphi(1))|=|4-3|=1.
$$
Only $\id$ and $n\mapsto1-n$ remain.
For $A=Y$, we have $d_Y(0,2)=6$, whereas $d_Y(j,j+2)=5$ for $j\ne0$.
The bijection $\varphi$ in~\eqref{eq:maps} preserves $|n-m|$, so it must preserve $\{0,2\}$: otherwise
$$
 \dis R\ge|d_Y(0,2)-d_Y(\varphi(0),\varphi(2))|=|6-5|=1.
$$
Only $\id$ and $n\mapsto2-n$ remain.
Thus $\RR_{<1}(X,X)$ and $\RR_{<1}(Y,Y)$ are finite and consist of graphs of bijections.

Suppose that $M$ is a midpoint between $X$ and $Y$.
Applying Lemma~\ref{lem:local} with $c=1$ first to $A=X$ and then to $A=Y$, we obtain optimal correspondences $R\in\RR(X,M)$ and $T^{-1}\in\RR(Y,M)$, since $1/8<1/4$.
Their distortions are $\dis R=2\dGH{}(X,M)=1/4$ and $\dis T^{-1}=2\dGH{}(Y,M)=1/4$.
The inverse correspondence $T\in\RR(M,Y)$ has the same distortion as $T^{-1}$.
By the composition estimate~\eqref{eq:composition},
$$
 \dis(T\circ R)\le\frac14+\frac14=\frac12,
$$
contrary to the nonexistence of a correspondence between $X$ and $Y$ of distortion at most $1/2$.

Since every shortest curve passes through a midpoint between its endpoints, there is no shortest curve between $X$ and $Y$.
\end{proof}

\printbibliography
\end{document}